\documentclass{ifacconf}
\usepackage{xcolor}

\usepackage{natbib}        

\usepackage{amssymb}
\usepackage{eqnarray} 

\usepackage{amsmath}
\usepackage{graphicx}	
\usepackage{datetime2}
\usepackage{pdfpages}		
\usepackage{dsfont}   		
\usepackage{algorithm}
\usepackage{algpseudocode}

\usepackage{bbm}  

\newtheorem{Thm}{Theorem}

\newtheorem{Lem}{Lemma}
\newtheorem{Rmk}{Remark}
\newtheorem{Def}{Definition}
\newtheorem{Property}{Property}

\newtheorem{Prob}{Problem}   
 
\newenvironment{proof}{\textit{Proof:}}{\hfill$\square$}

\newcommand{\Nh}{n}

\newcommand{\R}{\mathbb{R}}
\newcommand{\E}{\mathbb{E}}

\newcommand{\KL}{\text{KL}}
\newcommand{\Ipdf}{g}

\newcommand{\constrEqSymb}{\mathbf{E}}
\newcommand{\constrEq}[2]{\constrEqSymb_{#1,#2}}

\begin{document}
\begin{frontmatter}

\title{On the synthesis of control policies from noisy example datasets: a probabilistic approach} 


\author[First]{Davide Gagliardi} 
\author[Second]{Giovanni Russo} 

\address[First]{School of Electrical and Electronic Engineering, University College Dublin, Ireland (e-mail: davide.gagliardi@ucd.ie).}
\address[Second]{Department of Information and Electronic Engineering and Applied Mathematics, University of Salerno, Italy (e-mail: giovarusso@ucd.ie).}

\begin{abstract}
In this note we consider the problem of synthesizing optimal control policies for a system from noisy datasets. We present a novel algorithm that takes as input the available dataset and, based on these inputs, computes an optimal policy for possibly stochastic and non-linear systems that also satisfies actuation constraints. The algorithm relies on solid theoretical foundations, which have their key roots into a probabilistic interpretation of dynamical systems. The effectiveness of our approach is illustrated by considering an autonomous car use case. For such use case, we make use of our algorithm to synthesize a control policy from noisy data allowing the car to merge onto an intersection, while satisfying additional constraints on the variance of the car speed.
\end{abstract}
 
\end{frontmatter}

\section{Introduction}

A framework that is becoming particularly appealing to design control algorithms is that of devising the control policy from examples (or demonstrations), see e.g. \cite{Han_lIU_zHU_Pas_19,8550288} and references therein. At their roots these {\em control from demonstration techniques}, which are gaining considerable attention under the label of Inverse Reinforcement Learning (IRL), rely on Inverse Optimal Control and Optimization \cite{506395}. Today, IRL/control is recognized as an appealing framework to learn policies from {\em success stories} \cite{ARGALL2009469} and potential applications include planning \cite{doi:10.1177/0278364917745980} and preferences/prescriptions learning \cite{8796280}.

There is then no surprise that, over the years, a number of techniques have been developed to address the  problem of devising control policies from demonstrations, mainly in the context of Markov Decision Processes (MDPs) \cite{Sutton:1998:IRL:551283}. Results include \cite{Ratliff2009}, which leverages a linear programming approach, \cite{Ratliff:2006:MMP:1143844.1143936} which relies on a maximum margin approach, \cite{ziebart2008maximum} that makes use of the maximum entropy principle and \cite{Ramachandran:2007:BIR:1625275.1625692} that formalizes the problem via Bayesian statistics.
 
In this context, the main contributions of this extended abstract can be summarized as follows. First, we introduce an approach to synthesize control policies from examples which is based on the Fully Probabilistic Design (FPD)  \cite{Karny_M_Automatica_1996_Towards_Fully_Prob,Karny_M+Guy_TV_Sys&Ctr_lett_2006,Herzallah_R_JNeurNet_2015,Pegueroles_G+Russo_G_ECC19_confid,KARNY2012105}. This approach formalizes the control problem as an optimization problem where the Kullback-Leibler Divergence (see Section \ref{sec:KL_divergence}) between an {\em ideal} probability density function (pdf, obtained from e.g. demonstrations) and the pdf modeling the system/plant is minimized. The main technical novelty of our results with respect to the classic works on FPD lies in the fact that we explicitly embed actuation constraints in our formulation, thus solving an optimization problem where the Kullback-Leibler Divergence is minimized subject to constraints on the control variable. By relying on the FPD, one of the main advantages of our results over classic IRL/Control approaches is that policies can be synthesized from noisy data without requiring any assumption on the linearity of the system. The system can in fact be a general stochastic nonlinear dynamical system. Moreover, by embedding actuation constraints into the problem formulation and by solving the resulting optimization, we can export the policy that has been learned on other systems that have different actuation capabilities. As an additional contribution, we devise from our theoretical results an algorithmic procedure. The key reference applications over which the algorithm was tested involved an autonomous driving use case and full results are presented here.

\section{Mathematical Preliminaries}  
\subsection{Notation}
\label{ss:Notat}

Sets, as well as operators, are denoted by {\em calligraphic} characters, while vector quantities are denoted in {\bf bold}. Let $n_z$ be a positive integer and consider the measurable space $(\mathcal{Z},\mathcal{F}_z)$, with $\mathcal{Z}\subseteq\R^{n_z}$ and with $\mathcal{F}_z$ being a $\sigma$-algebra on $\mathcal{Z}$. 
Then, the random vector (i.e. a multidimensional random variable) on $(\mathcal{Z},\mathcal{F}_z)$ is denoted by $\mathbf{Z}$ and its realization is denoted by $\mathbf{z}$ (in the paper, we use the convention that these random vectors are row vectors). The \textit{probability density function} (or simply \textit{pdf} in what follows) of a continuous $\mathbf{Z}$ is denoted by 
$f_\mathbf{Z}(\mathbf{z})$. For notational convenience, whenever it is clear from the context, we omit the argument and/or the subscript of the pdf. Hence, the support of $f:=f_\mathbf{Z}(\mathbf{z})$ is denoted by $\text{S}\left(f\right)$ and, analogously, the  expectation of a function $\mathbf{h}(\cdot)$ of $\mathbf{Z}$ is indicated with $\E_{{f}}[\mathbf{h}(\mathbf{Z})]$ 
ad defined as $\E_{{f}}[\mathbf{h}(\mathbf{Z})] := \int_{\text{S}\left(f\right)}\mathbf{h}(\mathbf{z})f(\mathbf{z})d\mathbf{z}$. We also remark here that whenever we apply the averaging operator to a given function, we use an upper-case letter for the function argument as this is a random vector. The \textit{joint} pdf of two random vectors, say $\mathbf{Z}$ and $\mathbf{Y}$,is denoted by  $f_\mathbf{[Z,Y]}(\mathbf{z},\mathbf{y})$ and abbreviated with $f(\mathbf{z,y})$.
The \textit{conditional} probability density function ( or {\em cpdf} in what follows) of $\mathbf{Z}$ with respect to the random vector $\mathbf{Y}$ is denoted by $f\left( \mathbf{z}| \mathbf{y} \right)$ and, whenever the context is clear, we use the shorthand notation $\tilde{f}_\mathbf{Z}$. Finally, given $\mathcal{Z}\subseteq\R^{n_z}$, its \textit{indicator function} is denoted by $\mathds{1}_{\mathcal{Z}}(\mathbf{z})$. That is, $\mathds{1}_{\mathcal{Z}}(\mathbf{z}) =1$, $\forall \mathbf{z} \in \mathcal{Z}$ and $0$ otherwise. We also make use of the internal product between tensors, which is denoted by $\langle  \cdot , \cdot \rangle$.

\subsection{The Kullback-Leibler divergence}\label{sec:KL_divergence}
The control problem considered in this paper will be stated (see Section \ref{sec:control_problem}) in terms of the Kullback-Leibler (KL, \cite{KL_51}) divergence, formalized with the following:

\begin{Def}[Kullback-Leibler(KL) divergence]\label{def:KLdiv}	
Consider two pdfs, $\phi:=\phi_{\mathbf{Z}}(\mathbf{z})$ and $g:=g_{\mathbf{Z}}(\mathbf{z})$, with $ \phi$ being absolutely continuous with respect to $g$. Then, the \KL-divergence of $\phi$ with respect to $g$ is 
	\begin{equation}
	\label{eq:DKL_def}
	\mathcal{D}_{\KL}
	\left(\phi ||g \right):
	= \int_{\text{S}(\phi)} \phi \; \ln\left( \frac{\phi}{g}\right)\,d\mathbf{z}.
	\end{equation}
\end{Def}

Intuitively, $\mathcal{D}_{\KL}\left(\phi || g \right)$ is a measure of how well $\phi$ approximates $g$.  We now give give a property of the KL-divergence, the \KL-divergence splitting property, which is used in the proof of Theorem\ref{theo:ctrlConstr}.
\begin{Property}\label{proper:KLsplit} Let $\phi$ and $g$ be two pdfs of the random vector $[\mathbf{Z},\mathbf{Y}]$, with $\mathbf{Z}$ and $\mathbf{Y}$ being random vectors of dimensions $n^Z$ and $n^Y$, respectively. Then, the following {\em splitting} rule holds:
\begin{equation}
	\label{eq:KLgenSplit}
	\begin{array}{l}
		\mathcal{D}_{\KL}	
		\left(
		\phi (\mathbf{y},\mathbf{z}) || \, g (\mathbf{y},\mathbf{z}) 
		\right) =	\\
		\quad 
		\mathcal{D}_{\KL}
		\left( \phi (\mathbf{y}) || \, g (\mathbf{y} ) \right) + 
		\mathbb{E}_{ \phi ( \mathbf{Y} ) }
		\left[	
		\mathcal{D}_{\KL} \left( 	\phi(\mathbf{z}|\mathbf{Y}) || \, g(\mathbf{z}|\mathbf{Y}) \right)
		\right]	
	\end{array}		
\end{equation}
\proof 
The proof follows from the definition of $\mathcal{D}_{\KL}$, the \textit{conditioning} and \textit{independence} rules for pdfs. A self-contained proof of this technical result is reported in the appendix. 
\qed 	
\end{Property}

\section{Formulation of the Control Problem}
\label{ss:Assumpt}

Let: (i) $\mathcal{K} :=\lbrace k \rbrace_{k=1}^n$, $\mathcal{K}_0 := \mathcal{K} \cup \lbrace 0 \rbrace$ and $\mathcal{T}:= \lbrace t_k:  k \in\mathcal{K}_0\rbrace$ be the time horizon over which the system is observed; (ii) $\mathbf{x}_k\in\R^{d_x}$  and $\mathbf{u}_k\in\R^{d_u}$  be, respectively, the system state and input at time $t_k\in\mathcal{T}$; (ii) $\mathbf{d}_k:=(\mathbf{x}_k,\mathbf{u}_k)$ be the data collected from the system at time $t_k\in\mathcal{T}$ and  $\mathbf{d}^k$ the data collected from $t_0\in\mathcal{T}$ up to time $t_k\in\mathcal{T}$ ($t_k>t_0$). As shown in e.g. \cite{Peterka_V_Bayesian_Approach_to_sys_ident_1981}, the system behavior can be described via the joint pdf of the observed data, say $f(\mathbf{d}^n)$. Then, as shown in the same paper, the application of the chain rule for probability density functions leads to the following factorization for $f(\mathbf{d}^n)$:
\begin{equation}
\label{eq:CL_red_v1}
f\left( \mathbf{d}^n \right)  
=
\prod_{k\in \mathcal{K}} 
f	\left(	\mathbf{x}_k | \mathbf{u}_k, \mathbf{x}_{k-1} \right)
f	\left( \mathbf{u}_k | \mathbf{x}_{k-1} \right) 
f 	\left(	\mathbf{x}_0 	\right).
\end{equation}
Throughout this work we refer to (\ref{eq:CL_red_v1}) as the \textit{probabilistic description of the closed loop system}, or we simply say that (\ref{eq:CL_red_v1}) is our\textit{ closed loop} system. 

\begin{Rmk}
The cpdf $f\left(\mathbf{x}_k | \mathbf{u}_k, \mathbf{x}_{k-1} \right)$  describes the system behavior at time $t_k$, given the previous state and the input at time $t_k$. In turn, the input is also generated from the cdpf $f	\left( \mathbf{u}_k | \mathbf{x}_{k-1} \right)$, which is a {\em randomized control policy}, returning the input given the previous state. Finally, we also note that the initial conditions are embedded in the probabilistic system description through the prior $f \left(	\mathbf{x}_0 \right)$.
\end{Rmk}

In the rest of the paper we use the following \textit{shorthand} notations: $\tilde f_{\mathbf{X}}^k:=f	\left(	\mathbf{x}_k | \mathbf{u}_k, \mathbf{x}_{k-1} \right)$, 
$\tilde f_{\mathbf{U}}^k:=f	\left( \mathbf{u}_k | \mathbf{x}_{k-1} \right)$, $f_0:=f 	\left(	\mathbf{x}_0\right)$ and $f^n:= f\left( \mathbf{d}^n \right)$. Hence, (\ref{eq:CL_red_v1}) can be compactly written as
\begin{equation}
\label{eq:CL_red_v1_short}
f^n
=
\prod_{k\in \mathcal{K}} \tilde f_{\mathbf{X}}^k\tilde f_{\mathbf{U}}^kf_0 = \tilde f^n f_0, \ \ \ \tilde f^n:=\prod_{k\in \mathcal{K}} \tilde f_{\mathbf{X}}^k\tilde f_{\mathbf{U}}^k.
\end{equation}

\subsection{The control problem}\label{sec:control_problem}
Our goal is to synthesize, from an example dataset, say  $\mathbf{d}_e^n $, the control pdf $f\left( \mathbf{u}_k | \mathbf{x}_{k-1} \right)$ that allows the closed-loop system (\ref{eq:CL_red_v1_short})  to {\em achieve} the demonstrated behavior, subject to its actuation constraints. As in \cite{Karny_M_Automatica_1996_Towards_Fully_Prob,Quinn_A_Karny_M_Guy_FPD_of_Hierarchical_Bayes_mdl_InfoSci_2016,Pegueroles_G+Russo_G_ECC19_confid,Karny_M+Guy_TV_Sys&Ctr_lett_2006,Herzallah_R_JNeurNet_2015} the behavior illustrated in the example dataset can be specified through the reference pdf $\Ipdf \left( \mathbf{d}_e^n \right) $ extracted from the example dataset (as e.g. its empirical distribution). Following the chain rule for pdfs we have $\Ipdf \left( \mathbf{d}_e^n \right)  
	:=
	\prod_{k\in \mathcal{K}} 
	\Ipdf\left(	\mathbf{x}_k | \mathbf{u}_k, \mathbf{x}_{k-1} \right)
	\Ipdf\left( \mathbf{u}_k | \mathbf{x}_{k-1} \right) 
	\Ipdf\left(	\mathbf{x}_0\right)$. Again, by setting 
	$\tilde g_{\mathbf{X}}^k:=g	\left(	\mathbf{x}_k | \mathbf{u}_k, \mathbf{x}_{k-1} \right)$, 
	$\tilde g_{\mathbf{U}}^k:=g	\left( \mathbf{u}_k | \mathbf{x}_{k-1} \right)$, $g_0:=g 	\left(	\mathbf{x}_0\right)$
 and $g^n:=g\left( \mathbf{d}_e^n \right)$ we get:
\begin{equation}
\label{eq:CL_red_ref_short}
g^n
=
\prod_{k\in \mathcal{K}} \tilde g_{\mathbf{X}}^k\tilde g_{\mathbf{U}}^kg_0 = \tilde g^n g_0,
\end{equation}
where $\tilde g^n:=\prod_{k\in \mathcal{K}} \tilde g_{\mathbf{X}}^k\tilde g_{\mathbf{U}}^k$.

The control problem can then be recast as the problem of designing $f\left( \mathbf{u}_k | \mathbf{x}_{k-1} \right)$ so that $f^n$ approximates $g^n$. This leads to the following formalization:

\begin{Prob}\label{prob:Main_Constr_Ctrl}
Determine the sequence of cpdfs, say
$\left \lbrace 
	\left( \tilde{f}_\mathbf{U}^k \right)^* 
	\right\rbrace_{k \in \mathcal{K}}$,  
solving the nonlinear program
\begin{equation}\label{eq:ProbMain}	
\begin{aligned}
	& \underset{{\left \lbrace\tilde{f}_\mathbf{U}^k\right \rbrace_{k \in \mathcal{K}}}}{\min}	
	& & \mathcal{D}_{\KL}\left(f^n||\Ipdf^n\right)\\
	& \text{s.t.}
	& & \E_{\tilde{f}_\mathbf{U}^k}
	\left[
	\tilde{\mathbf{h}}_{\mathbf{u},k} \left(\mathbf{U}\right)
	\right] =  \tilde{\mathbf{H}}_{\mathbf{u},k}, 
		\quad k \in \mathcal{K} ,
\end{aligned}
\end{equation}
where the constraints are algebraically independent.
\end{Prob}
 In Problem \ref{prob:Main_Constr_Ctrl}, the constraints are formalized as expectations. We note that these constraints can be equivalently written as
$\int_{\text{S}(\tilde{f}_\mathbf{U}^k )} \tilde{f}_\mathbf{u}^k \; \tilde{\mathbf{h}}_{\mathbf{u},k}\left( \mathbf{u} \right)\; d\mathbf{u} =  \tilde{\mathbf{H}}_{\mathbf{u},k}$. Also, the constraints of the program are time-varying and the number of constraints can change over time (the number of constraints at time $t_k$ is denoted by $c_{\mathbf{u},k}$). Indeed, in the constraints of (\ref{eq:ProbMain}): (i) $\tilde{\mathbf{H}}_{\mathbf{u},k}$ is a (column) vector of coefficients, i.e. $\tilde{\mathbf{H}}_{\mathbf{u},k}:= \left[H_{\mathbf{u},0,k}, \mathbf{H}_{\mathbf{u},k}^T\right]^T$ and $\tilde{\mathbf{h}}_{\mathbf{u},k} \left( \mathbf{z} \right)  :=\left[h_{\mathbf{u},0,k}, \mathbf{h}_{\mathbf{u},k}^T\right]^T(\mathbf{z})$; (ii) $\mathbf{H}_{\mathbf{u},k} \in \mathbb{R}^{c_{\mathbf{u},k}}$ and $\mathbf{h}_{\mathbf{u},k}: \text{S}(\tilde{f}_\mathbf{u}^k ) \mapsto \mathbb{R}^{c_{\mathbf{u},k}}$; (iii) $H_{\mathbf{u},0,k} :=1$ and $h_{\mathbf{u},0,k} \left( \mathbf{z} \right) := \mathds{1}_{\mathcal{U}_k}  \left( \mathbf{z} \right)$ ensure that the solution of the program is a cpdf. Finally, in  Problem \ref{prob:Main_Constr_Ctrl} we assume that the constraints are algebraically independent. The notion of algebraically independent constraints is formalized next.
	\begin{Def}
		\label{def:algebrIndep}
Let $\mathbf{Z}$ be a random vector with underlying pdf $f_\mathbf{Z} \left( \mathbf{z} \right)$ and support $\mathcal{Z}$. 
		A set of functions $\mathbf{h} : \mathcal{Z} \mapsto \mathbb{R}^{c_{\mathbf{z}}}$ is said to be \emph{algebraically independent} if there exists a subset, say $S \subset \mathcal{Z}$, with non-zero measure (i.e. $\int_S\, d\mathbf{z}>0$) and such that:
		\begin{equation}\label{eqn:algebraically_independent}
		\exists	S \subset \mathcal{Z}:
		\quad
		\int_S 
		\langle 
		\mathbf{v},\,\mathbf{h}\left( \mathbf{z} \right) 
		\rangle^2
		\, d\mathbf{z}
		>0, 
		\quad 
		\forall 
		\mathbf{v}
		\in \mathbb{R}^{c_{\mathbf{z}}}\backslash{\mathbf{0}}	
		\end{equation}	
	\end{Def}
In what follows, we simply say that a set of equations (or constraints) of the form of (\ref{eqn:algebraically_independent}) is algebraically independent if the above definition is satisfied. As shown in \cite{Guilleminot_J_Soize_C_JE_2013_OnStatDependenceOfMaterilSymmetry}, the assumption that the contraints are algebraically independent ensures that Problem \ref{prob:Main_Constr_Ctrl} is well posed.

\section{Technical Results}

We now introduce the main technical results of this paper. The key result behind the algorithm of  Section \ref{sec:algorithm} is Theorem \ref{theo:ctrlConstr}. The proof of this result, given in this section, makes use of three technical lemmas (i.e. Lemma \ref{lem:Constrained_KL}, Lemma \ref{lem:LM_optim} and Lemma\ref{proper:KLdiv_split_nD_fn}). 

\begin{Lem}\label{lem:Constrained_KL}
	
Let: (i) $\mathbf{Z}$ be a random 
	vector on the measurable space $\left(\mathcal{Z},\mathcal{F}_z\right)$; (ii) $f := f_\mathbf{Z} (\mathbf{z})$, $g:= g_\mathbf{Z} (\mathbf{z})$ be two probability distributions  over $\left(\mathcal{Z},\mathcal{F}_z\right)$; (iii) $\alpha: \mathcal{Z} \mapsto \mathbb{R}_0^+$ be a \textit{nonnegative} function of $\mathbf{Z}$, \textit{integrable} under the measure given by $ f_\mathbf{Z}(\mathbf{z})$. Assume that  $f_\mathbf{Z}(\mathbf{z})$ satisfies the following set of algebraically independent equations:
	\begin{equation}
	\label{eq:gen_constr}
	\begin{array}{lrl}
	\int{f_\mathbf{Z} \left( \mathbf{z} \right) \;	
		\tilde{\mathbf{h}}\left( \mathbf{z} \right)\; d\mathbf{z}} 
	= 
	\tilde{\mathbf{H}},
	\end{array} 	
	\end{equation}	
	where: (i) $\tilde{\mathbf{h}}(\mathbf{z}):= \left[h_0, \mathbf{h}^T \right]^T(\mathbf{z})$, with $h_0 (\mathbf{z}) :=  \mathds{1}_{\mathcal{S}  \left( \mathbf{Z} \right) } (\mathbf{z})$ and $\mathbf{h}:\mathcal{Z} \mapsto \mathbb{R}^{c_{\mathbf{z}}}$ being a measurable map; (ii) $\tilde{\mathbf{H}}(\mathbf{z}) :=\left[H_0, \mathbf{H}^T\right]^T$ with $H_0:= 1$ and  $\mathbf{H} \in \mathbb{R}^{c_{\mathbf{z}}}$ being a vector of constants. Then:
\begin{enumerate}
\item the solution of the constrained optimization problem
	\begin{equation}\label{eqn:probl_Lemma_1}
	\begin{aligned}
		& \underset{f_\mathbf{Z}}{\min}	
		& &\mathcal{L}(f) 
		& \text{s.t.}
		& & \text{constraints in \eqref{eq:gen_constr}} 
	\end{aligned}
	\end{equation}
	with 
	\begin{equation}
	\mathcal{L}(f) := \mathcal{D}_{\KL} \left(f || g \right)	+ 
	\int f_\mathbf{Z} (\mathbf{z}) \; \alpha \left( \mathbf{z} \right) \, d\mathbf{z} 
	\end{equation}	
	is the pdf
	\begin{equation}
	\label{eq:Lem1_opt_fu_sol}
	f^\ast:=f^\ast_\mathbf{Z}\left( \mathbf{z} \right) 
	= 
	\frac{ 
		g\left( \mathbf{z} \right)  
		e^{- 
			\lbrace
			\alpha \left( \mathbf{z} \right) +
			\langle
			\boldsymbol{\lambda}^*,
			\mathbf{h}\left( \mathbf{z} \right)
			\rangle
			\rbrace 
		}
	}{e^{1+\lambda^*_0}}.
	\end{equation}	
	In \eqref{eq:Lem1_opt_fu_sol} $\lambda^\ast_0$ and $\boldsymbol{\lambda}^\ast=[\lambda^\ast_1,\dots,\lambda^\ast_{c_{\mathbf{z}}}]^T$ are the \textit{Lagrange multipliers} associated to the constraints;
\item moreover, the corresponding minimum is:
\begin{equation}
\label{eq:minVal}
\mathcal{L}^\ast :
= \mathcal{L} \left( f^* \right)  
=-\left( 1+\lambda^\ast_0 +	\langle \boldsymbol{\lambda}^\ast, \mathbf{H}\rangle	\right).
\end{equation}	
\end{enumerate}
\end{Lem}
\proof $\;$ See the appendix. \qed \\
 
Note that, in Lemma \ref{lem:Constrained_KL}, the optimal solution $f_\mathbf{Z}^\ast(\mathbf{z})$ depends on the Lagrange multipliers (LMs) $\lambda_0^\ast$ and $\boldsymbol{\lambda}^\ast$. The first LM, i.e. $\lambda_0^\ast$, can be obtained by integration, i.e. by imposing that $e^{1+\lambda_0^\ast}$ normalizes $f_\mathbf{Z}^\ast(\mathbf{z})$ in (\ref{eq:Lem1_opt_fu_sol}). With the next result, we propose a strategy for finding the LMs $\boldsymbol{\lambda}^\ast$. 
In particular, the idea is to recast the problem of finding the solutions of non-linear equations as a minimization problem. 
In general, the approach can be also used to fit the parameters of a pdf so that it meets a set of pre-specified constrains (for example, to find pdfs that satisfy the Maximum Entropy principle \cite{Guilleminot_J_Soize_C_JE_2013_OnStatDependenceOfMaterilSymmetry}).

\begin{Lem}
\label{lem:LM_optim}
Let: (i) $\mathcal{Z}\subseteq \mathbb{R}^{n_\mathbf{z}}$ and $\tilde{\Theta} \subseteq \mathbb{R}^{n_\mathbf{z}}$; 
	(ii) $\hat{f}_1:\mathcal{Z}\mapsto \hat{f}_1\left( \mathbf{z} \right)$ be a positive and integrable function on $\mathcal{Z}$;
	(iii) $\hat{f}_2: (\mathcal{Z}\times \tilde{\Theta})\mapsto 
	\hat{f}_1 \left( \mathbf{z} \right)	e^{- 
		\langle \tilde{\boldsymbol{\theta}}  
		,\; 	\tilde{\mathbf{h}} \left( \mathbf{z} \right)
		\rangle
		}$, where $\tilde{\mathbf{h}}=\left[\tilde{\mathbf{h}}_1\left( \mathbf{z} \right),\ldots,\tilde{\mathbf{h}}_{c_{\mathbf{z}}}\left( \mathbf{z} \right)\right]^T : \mathcal{Z} \mapsto \mathbb{R}^{c_{\mathbf{z}}}$ are algebraically independent functions. Consider the constraints defined by the set of the following equations: 
	\begin{equation}
	\label{eq:constr_tilde_Theta}
	\int_{\mathcal{Z}}	 
	\hat{f}_2 
	\left( 
	\mathbf{z},  \tilde{\boldsymbol{\theta}}  
	\right) 	
	\tilde{\mathbf{h}}_i\left( \mathbf{z} \right)\; d\mathbf{z} 
	= 
	\tilde{\mathbf{H}}_i ,\;
	\quad 
	i = 1,\dots,c_{\mathbf{z}},
	\end{equation}
	where $\tilde{\mathbf{H}}:=\left[\tilde{\mathbf{H}}_1,\ldots,\tilde{\mathbf{H}}_{c_{\mathbf{z}}}\right]^T\in \mathbb{R}^{c_{\mathbf{z}}}$. Then, the unique solution, say  $\tilde{\boldsymbol{\theta}}^\ast$, of the minimization problem
\begin{equation}
\label{eq:Index_Inique}
\underset{\tilde{\boldsymbol{\theta}}}{\min}
\;\mathcal{J}
\left( \tilde{\boldsymbol{\theta}} \right),
\end{equation}
with 
$\mathcal{J} \left( \tilde{\boldsymbol{\theta}} \right):
=
\langle 
\tilde{\boldsymbol{\theta}}  
,\, 	\tilde{\mathbf{H}}
\rangle
+
\int_{\mathcal{Z}} 
\hat{f}_2 \left( \mathbf{z},  \tilde{\boldsymbol{\theta}}  \right) 	\, d\mathbf{z} 
$
is also a solution of \eqref{eq:constr_tilde_Theta}. 
\end{Lem}
\proof  $\;$ See the appendix \qed \\

Finally,  we introduce here the following technical lemma that is used in the proof of Theorem \ref{theo:ctrlConstr}.
\begin{Lem}\label{proper:KLdiv_split_nD_fn}	
Let $f^n$ and $g^n$ be the pdfs defined in 
\eqref{eq:CL_red_v1} and \eqref{eq:CL_red_ref_short}, respectively. Then:
\begin{equation}
\label{eq:KLdiv_split_nD_fn}
\mathcal{D}_{\KL}
\left( f^n || \Ipdf^n \right)  		
=
\mathcal{D}_{\KL}
\left( f^{n-1} ||\Ipdf^{n-1} \right)  
+
\mathbb{E}_{f^{n-1}}
\left[	
\mathcal{D}_{\KL} 
\left(
\tilde{f}^{n}|| \tilde{\Ipdf}^{n}
\right)
\right]	
\end{equation}	
\proof 
The result is obtained from Property~\ref{proper:KLsplit} (see the appendix for a proof of this property) by setting
$\mathbf{Y} := [\mathbf{X}_0,\mathbf{U}_1,\mathbf{X}_1,\dots, \mathbf{U}_{n-1}, \mathbf{X}_{n-1} ]$ 
and
$\mathbf{Z}:=[\mathbf{U}_{n}, \mathbf{X}_{n} ]$
\qed 
\end{Lem}

The main result behind the algorithm of Section \ref{sec:algorithm}, the proof of which makes use of the above technical results, is presented next.
\begin{Thm}\label{theo:ctrlConstr}
The solution,
$\left( \tilde{f}_\mathbf{U}^k \right)^\ast = f^\ast \left( \mathbf{u}_k | \mathbf{x}_{k-1} \right)$, of the control Problem \ref{prob:Main_Constr_Ctrl} is
\begin{equation}
\label{eq:opt_ctrl}
\left( \tilde{f}_\mathbf{U}^k \right)^* = 
\tilde{\Ipdf}_\mathbf{U}^k 
\frac{e^{-\lbrace 
	\hat{\omega} \left( \mathbf{u}_k,\,\mathbf{x}_{k-1} \right) +
	\langle
	\boldsymbol{\lambda}^*_{\mathbf{u},k}, \mathbf{h}_{\mathbf{u},k}\left( \mathbf{u}_k \right)
	\rangle
	\rbrace}}{e^{1+\lambda^*_{\mathbf{u},0,k}}},
\end{equation}
where:
\begin{enumerate}
\item $\hat{\omega}(\cdot,\cdot)$ is generated via the backward recursion
	\begin{equation}
		\label{eq:omega_defdef}
			\hat{\omega} \left( \mathbf{u}_k,\,\mathbf{x}_{k-1} \right)  = 	
			\hat{\alpha} \left( \mathbf{u}_k,\,\mathbf{x}_{k-1} \right) +
			\hat{\beta} \left( \mathbf{u}_k,\,\mathbf{x}_{k-1} \right), 
		\end{equation}
with
\begin{equation}
\label{eq:def:alpha_beta}
\begin{split}
\hat{\alpha} \left( \mathbf{u}_k,\,\mathbf{x}_{k-1} \right)
	& :=
	\mathcal{D}_{\KL}
			\left(
			\tilde{f}^k_{\mathbf{X}}
			|| 
			\tilde{\Ipdf}^k_{\mathbf{X}}
			\right)\\
	\hat{\beta} \left( \mathbf{u}_k,\,\mathbf{x}_{k-1} \right)
	& := 
			- \mathbb{E}_{\tilde{f}^k_\mathbf{X}}
			\left[	
			\ln 
			\hat{\gamma}
			\left(
			\mathbf{X}_{k}
			\right)
			\right],
		\end{split}
	\end{equation}
with terminal conditions $\hat{\beta} \left( \mathbf{u}_n,\,\mathbf{x}_{n-1} \right) = 0$ and $\hat{\alpha} \left( \mathbf{u}_n,\,\mathbf{x}_{n-1} \right) = \mathcal{D}_{\KL}\left(\tilde{f}^n_{\mathbf{X}}||\tilde{\Ipdf}^n_{\mathbf{X}}
\right)$;	
\item $\hat{\gamma}\left( \cdot \right)$ in \eqref{eq:def:alpha_beta} is given by
\begin{equation}
	\label{eq:def_gamma}
	\ln 
	\hat{\gamma}
	\left(
	\mathbf{x}_{k-1}
	\right)
	:=	\left[
	\sum_{i=0}^{c_{\mathbf{u}}}
	\ln
	\left(
	\hat{\gamma}_{\mathbf{u},i,k}				
	\left(
	\mathbf{x}_{k-1}
	\right)
	\right) 	
	\right], 
\end{equation}
with
		\begin{equation}
		\label{eq:def_gamma_0}
		\hat{\gamma}_{\mathbf{u},0,k} \left( \mathbf{x}_{k-1} \right)
		= \exp{ \lbrace	\lambda^*_{\mathbf{u},0,k} + 1 \rbrace },
		\end{equation}
		and
		\begin{equation}
		\label{eq:def_gamma_i}
		\hat{\gamma}_{\mathbf{u},i,k}				
		\left(
		\mathbf{x}_{k-1}
		\right) :=  
		\exp{ \lbrace \lambda^*_{\mathbf{u},i,k} \mathbf{H}_{\mathbf{u},i,k}	\rbrace}  	\quad i = 1,\dots,c_{\mathbf{u}},
		\end{equation}	
		with terminal conditions
		$\hat{\gamma}_{\mathbf{u},0,n}   \left(\mathbf{x}_{n-1} \right) = 1$, i.e. $\lambda^*_{\mathbf{u},0,n} = 0$, and  $\lambda^\ast_{\mathbf{u},i,n} = 0$, $i = 1,\ldots,n$; 	
\item $\lambda^\ast_{\mathbf{u},0,k}$  and $\boldsymbol{\lambda}^\ast_{\mathbf{u},k} = \left[\lambda_{\mathbf{u},1,k}^\ast,\dots,\lambda_{\mathbf{u},c_{\mathbf{u}},k}^\ast \right]$ in (\ref{eq:opt_ctrl}) are the Lagrange multipliers (LMs) associated to the constraints at time $t_k$. In particular, 
\begin{equation*}
\begin{split}
		& \lambda^\ast_{\mathbf{u},0,k}
		 =\\
		 &
		\ln
		\left[ 
		\int \tilde{\Ipdf}_\mathbf{U}^k 
		\left(
e^{-\lbrace 
				\hat{\omega} \left( \mathbf{u}_k,\,\mathbf{x}_{k-1} \right) +
				\langle
				\boldsymbol{\lambda}^*_{\mathbf{u},k}, \mathbf{h}_{\mathbf{u},k}\left( \mathbf{u}_k \right)
				\rangle
				\rbrace}
		\right)	\; 
		d\mathbf{u}_k 
		\right] - 1,
		\end{split}
		\end{equation*}
while all the other LMs can be obtained numerically (via e.g. Lemma \ref{lem:LM_optim}).	
\end{enumerate}	
Moreover, the corresponding minimum at time $k$ is given by:
\begin{equation}
\label{eq:ctrl_contrst_minVal}
B^*_k 
:= 	
- \mathbb{E}_{p^{k-1}_\mathbf{X}}
\left[
\ln \hat{\gamma} \left( \mathbf{X}_{k-1} \right) 	
\right]. 	
\end{equation}
where $p_{\mathbf{X}}^{k}$ denotes the pdf of the state at time $t_k$ 
(i.e. $p_{\mathbf{X}}^{k}:= f\left(\mathbf{x}_{k}\right)$). 
	
\end{Thm}
\proof
For notational convenience, we use the shorthand notation $\left \lbrace
	\constrEq{\mathbf{u}}{k}
	\right \rbrace$ to denote the set of constraints of Problem \ref{prob:Main_Constr_Ctrl} at time $t_k$. We also denote by $\left \lbrace
	\constrEq{\mathbf{u}}{k}\right\rbrace_{\mathcal{K}}$ the set of constraints over the whole time horizon $\mathcal{K}$ and $\left \lbrace
	\constrEq{\mathbf{u}}{k}
	\right \rbrace_{k = 1}^{n-1}$ to denote the constraints from  $t_1$ up to time $t_{n-1}$.


Note that, following Lemma \ref{proper:KLdiv_split_nD_fn}, Problem \ref{prob:Main_Constr_Ctrl} can be re-written as follows:
\begin{equation}
\label{eq:SplitLast_red}
\begin{array}{l}
\underset{	
	\footnotesize
	\begin{array}{l}
	\left\lbrace
	\tilde{f}_\mathbf{U}^k
	\right\rbrace_{k \in \mathcal{K}} \\
	\text{ s.t.:} \\ 
	\left \lbrace
	\constrEq{\mathbf{u}}{k}
	\right \rbrace_{k \in \mathcal{K}}	
	\end{array}
}{\min} 
\mathcal{D}_{\KL} \left( f^{n} || \Ipdf^{n} \right) =
\\
\qquad
=	
\underset{
	\footnotesize
	\begin{array}{l}
	\left\lbrace
	\tilde{f}_\mathbf{U}^k
	\right\rbrace_{k = 1}^{n-1} \\
	\text{ s.t.:} \\ 
	\left \lbrace
	\constrEq{\mathbf{u}}{k}
	\right \rbrace_{k = 1}^{n-1}
	\end{array}
}{\min} 
\left\lbrace
\mathcal{D}_{\KL} \left( f^{n-1} || \Ipdf^{n-1} \right) 
+ B^*_n
\right\rbrace 	
\\
\end{array}
\end{equation}
where:
\begin{equation}
\label{eq:BnStar_Ctrl_constr}
B^*_n 
:=
\underset{
		\footnotesize
		\begin{array}{l}
		\tilde{f}_\mathbf{U}^n \\
		\text{ s.t.:} \\ 
		\constrEq{\mathbf{u}}{n}
		\end{array}
	}{\min} B_n
\, ,
\quad
B_n 
:=
\mathbb{E}_{f^{n-1}}
\left[
\mathcal{D}_{\KL}
\left(
\tilde{f}^n ||
\tilde{\Ipdf}^n
\right)
\right].
\end{equation}
That is, Problem \ref{prob:Main_Constr_Ctrl} can be approached by solving first the optimization of the last time-instant of the time-horizon $\mathcal{K}$ 
(the term $B_n$ in \eqref{eq:SplitLast_red}) and then by taking into account the result from this optimization problem in the optimization up to the instant $t_{n-1}$. Now we focus on the sub-problem:
\begin{equation}
\label{eq:Prob_Bn_StarCtrl_constr}
B^*_n : =
\underset{
	\footnotesize
	\begin{array}{l}
	\tilde{f}_\mathbf{U}^n \\
	\text{ s.t.:} \\ 
	\constrEq{\mathbf{u}}{n}
	\end{array}
}{\min} B_n
\end{equation}
For this problem, we first observe that the following equality is satisfied for the term $B_n$:
\begin{equation}
\label{eq:Bn_reducion}
B_n 
=
\mathbb{E}_{f^{n-1}}
\left[
\mathcal{D}_{\KL}
\left(
\tilde{f}^n ||
\tilde{\Ipdf}^n
\right)
\right]
=
\mathbb{E}_{p^{n-1}_\mathbf{X}}
\left[
\mathcal{D}_{\KL}
\left(
\tilde{f}^n ||
\tilde{\Ipdf}^n
\right)
\right].
\end{equation}
Such equality was obtained by noting that $\mathcal{D}_{\KL} \left( \tilde{f}^n || \tilde{\Ipdf}^n \right)$ is only a function of the previous state (see also \cite{Karny_M_Automatica_1996_Towards_Fully_Prob}) and, for notational convenience, we rename it as $\hat{A}\left( \cdot \right)$. Hence, $B_n$ becomes 
\begin{equation}
\label{eq:Bn_Ctrl_constr}
B_n 
=
\mathbb{E}_{p^{n-1}_\mathbf{X}}
\left[
\mathcal{D}_{\KL}
\left(
\tilde{f}^n ||
\tilde{\Ipdf}^n
\right)
\right]
=
\mathbb{E}_{p^{n-1}_\mathbf{X}}
\left[ \hat{A}\left(\mathbf{X}_{n-1}\right) \right].
\end{equation}
Now, note that
\begin{equation}
\begin{array}{ll}
B^*_n:
&=
\underset{
	\footnotesize
	\begin{array}{l}
	\tilde{f}_\mathbf{U}^n \\
	\text{ s.t.:} \\ 
	\constrEq{\mathbf{u}}{n}
	\end{array}
}{\min} 
B_n
= 
\underset{
	\footnotesize
	\begin{array}{l}
	\tilde{f}_\mathbf{U}^n \\
	\text{ s.t.:} \\ 
	\constrEq{\mathbf{u}}{n}
	\end{array}
}{\min}
\mathbb{E}_{p^{n-1}_\mathbf{X}}
\left[ \hat{A}\left(\mathbf{X}_{n-1}\right) \right]=\\
& =
\mathbb{E}_{p^{n-1}_\mathbf{X}}
\left[
\underset{
	\footnotesize
	\begin{array}{l}
	\tilde{f}_\mathbf{U}^n \\
	\text{ s.t.:} \\ 
	\constrEq{\mathbf{u}}{n}
	\end{array}
}{\min}
\hat{A}\left(\mathbf{X}_{n-1}\right) 
\right]
= 
\mathbb{E}_{p^{n-1}_\mathbf{X}}
\left[
A^*_n
\right],
\end{array}
\end{equation}
where the above expression was obtained by using the fact that the expectation operator is linear and the fact that independence of the decision variable 
(i.e. $\tilde{f}_\mathbf{U}^n$) is independent on the pdf over which the expectation is performed (i.e. $p^{n-1}_\mathbf{X}$). This implies that, once we solve the problem
\begin{equation}
\label{eq:AnStars_Ctrl_constr}
A^\ast_n 
:=
\underset{
		\footnotesize
			\begin{array}{l}
				\tilde{f}_\mathbf{u}^n \\
				\text{ s.t.:} \\ 
				\constrEq{\mathbf{u}}{n}
			\end{array}
		}{\min} \hat{A}(\mathbf{x}_{n-1})
\end{equation}
for any fixed $\mathbf{x}_{n-1}$, then $B^*_n$ can be obtained by averaging $A^\ast_n $ over $p_{\mathbf{X}}^{n-1}$. We now focus on solving problem \eqref{eq:AnStars_Ctrl_constr}. In doing so, we first note that, following \eqref{eq:Bn_Ctrl_constr}, $\hat{A}(\mathbf{x}_{n-1})$ can be re-written as follows:
\begin{subequations}
	\begin{equation}
	\label{eq:An_expl_of_alpha}
	\hat{A}\left(\mathbf{x}_{n-1}\right)
	=
	\int
	\tilde{f}^n_\mathbf{U}
	\left[
	\ln
	\left(
	\frac{
		\tilde{f}^n_\mathbf{U}
	}{
		\tilde{\Ipdf}^n_\mathbf{U}
	}
	\right)
	+
	\hat{\alpha}
	\left( \mathbf{u}_n,\mathbf{x}_{n-1} \right)
	\right]\,
	d\mathbf{u}_{n},
	\end{equation}	
	\begin{equation}
	\hat{\alpha}
	\left( \mathbf{u}_n,\mathbf{x}_{n-1} \right)
	:=
	\mathcal{D}_{\KL} 
	\left( \tilde{f}^n_\mathbf{X} || \tilde{\Ipdf}^n_\mathbf{X} \right).
	\end{equation}
\end{subequations}
In turn, \eqref{eq:An_expl_of_alpha} can be compactly written as:
\begin{equation}
\hat{A}(\mathbf{x}_{n-1}) 
= 
\mathcal{D}_{\KL} \left( \tilde{f}_u^n || \tilde{\Ipdf}_u^n \right) 
+\int{\tilde{f}_\mathbf{U}^n \,
	\hat{\alpha} \left( \mathbf{u}_n,\mathbf{x}_{n-1} \right) 	\,
	d\mathbf{u}_n},
\end{equation}
where we used the definition of \KL-divergence.
Hence, Lemma \ref{lem:Constrained_KL} can be used to solve the optimization problem in.\eqref{eq:AnStars_Ctrl_constr}. Indeed by applying Lemma~\ref{lem:Constrained_KL} with $\mathbf{Z}=\mathbf{U}_n$, 
$f=\tilde{f}_\mathbf{U}^n$, 
$g=\tilde{g}_\mathbf{U}^n$,
$\mathbf{h} = \mathbf{h}_{\mathbf{u},n} $,
$\mathbf{H} = \mathbf{H}_{\mathbf{u},n} $
we get the following solution to \eqref{eq:AnStars_Ctrl_constr}:
\begin{equation}
\label{eq:opt_ctrl_n}
\left( \tilde{f}_\mathbf{U}^n \right)^\ast =
\tilde{\Ipdf}_\mathbf{U}^n 
\frac{e^{-\lbrace 
		\hat{\alpha} \left( \mathbf{u}_n,\,\mathbf{x}_{n-1} \right) +
		\langle
		\boldsymbol{\lambda}^\ast_{\mathbf{u},n},
		\mathbf{h}_{\mathbf{u},n} \left( \mathbf{u}_n \right)
		\rangle 
		\rbrace}}{e^{1+\lambda^\ast_{\mathbf{u},0,n}}}.
\end{equation}
In the above pdf, $\lambda^\ast_{\mathbf{u},0,n}$ and $\boldsymbol{\lambda}^\ast_{\mathbf{u},n}$ are the LMs at the last time instant, $t_n$. The LM $\lambda^\ast_{\mathbf{u},0,n}$ can be obtained by imposing a normalization condition to  (\ref{eq:opt_ctrl_n}). That is, $\lambda^\ast_{\mathbf{u},0,n}$ can be  found by imposing that
\begin{equation}
	\begin{array}{ll}
		\exp\{\lambda^*_{\mathbf{u},0,n} + 1\} & = \int \tilde{\Ipdf}_\mathbf{U}^n e^{-\lbrace \hat{\alpha} \left( \mathbf{u}_n,\,\mathbf{x}_{n-1} \right) +
			\langle
			\boldsymbol{\lambda}^*_{\mathbf{u},n},
			\mathbf{h}_{\mathbf{u},n}\left( \mathbf{u}_n \right)
			\rangle
			\rbrace}	
		\; d\mathbf{u}_n  \\
		& = \hat{\gamma}_{\mathbf{u},0,n}	
		\left( \mathbf{x}_{n-1} \right).
	\end{array}
\end{equation}
Also, following Lemma \ref{lem:Constrained_KL}, the minimum of the problem is given by:
\begin{equation}
\hat{A}^\ast_n =
- 	\left(
1 +	
\lambda^\ast_{\mathbf{u},0,n} +
\langle
\boldsymbol{\lambda}^\ast_{\mathbf{u},n},
\mathbf{H}_{\mathbf{u},n}
\rangle
\right)
\end{equation}
or equivalently
\begin{equation}
\label{eq:opt_An}
\hat{A}^\ast_n
=
-
\left[
\sum_{i=0}^{c_{\mathbf{u}}}
\ln
\left(
\hat{\gamma}_{\mathbf{u},i,n}				
\left(
\mathbf{x}_{n-1}
\right)
\right) 	
\right] 
=
- \ln \hat{\gamma} 
	\left(
	\mathbf{x}_{n-1}
	\right)
\end{equation}
where we have used the definitions \eqref{eq:def_gamma_0} and \eqref{eq:def_gamma_i} for $\hat{\gamma}_{\mathbf{u},i,n}, i = 0,\dots c_{\mathbf{u}}$. Therefore, the corresponding minimum value for $B_n$ is:
\begin{equation}
\label{eq:ctrl_contrst_min_Bn}
B^\ast_n
=
- \mathbb{E}_{p^{n-1}_\mathbf{X}}
\left[
\ln \hat{\gamma} 
\left(
\mathbf{X}_{n-1}
\right)
\right].
\end{equation}	
Note now that the solution we found to the problem in (\ref{eq:Prob_Bn_StarCtrl_constr}) only depends on $\mathbf{X}_{n-1}$ and therefore the original problem \eqref{eq:SplitLast_red} can be split as
\begin{equation}
\label{eq:Split2Last_red_explicit}
\begin{array}{ll}
\underset{
	\footnotesize
	\begin{array}{l}
	\left\lbrace
	\tilde{f}_\mathbf{U}^k
	\right\rbrace_{k = 1}^{n-1} \\
	\text{ s.t.:} \\ 
	\left \lbrace
	\constrEq{\mathbf{u}}{k}
	\right \rbrace_{k = 1}^{n-1}
	\end{array}
}{\min} 
&
\left\lbrace
\mathcal{D}_{\KL} \left( f^{n-1} || \Ipdf^{n-1} \right) 
+ B^*_n
\right\rbrace 	= \\
\quad =
\underset{
	\footnotesize
	\begin{array}{l}
	\left\lbrace
	\tilde{f}_\mathbf{U}^k
	\right\rbrace_{k = 1}^{n-2} \\
	\text{ s.t.:} \\ 
	\left \lbrace
	\constrEq{\mathbf{u}}{k}
	\right \rbrace_{k = 1}^{n-2}
	\end{array}
}{\min} 
&
\left\lbrace
\mathcal{D}_{\KL} \left( f^{n-2} || \Ipdf^{n-2} \right) 
+ B^*_{n-1}
\right\rbrace \\
\\
\end{array}
\end{equation}
where:
\begin{subequations}
\begin{equation}
\label{eq:Bnm1_Star_Ctrl_constr}
B^*_{n-1} 
:=
\underset{
	\footnotesize
	\begin{array}{l}
	\tilde{f}_\mathbf{U}^{n-1}\\
	\text{ s.t.:} \\ 
	\constrEq{\mathbf{u}}{n-1}
	\end{array}
}{\min} B_{n-1} , 
\end{equation}
\begin{equation}
\label{eq:Bnm1_Ctrl_constr}
B_{n-1} 
:=
\mathbb{E}_{f^{n-2}}
\left[
\mathcal{D}_{\KL}
\left(
\tilde{f}^{n-1} ||
\tilde{\Ipdf}^{n-1}
\right)
\right]+ B^\ast_n
\end{equation}
\end{subequations}
We approach the above problem in the same way we used to solve the problem in (\ref{eq:Prob_Bn_StarCtrl_constr}). The idea is now to find a function, $\hat{A}\left(\mathbf{x}_{n-2}\right)$, such that 
\begin{equation}
\label{eq:Bnm1_of_Anm1}
B_{n-1} 
=
\mathbb{E}_{p^{n-2}_\mathbf{X}}
\left[ \hat{A}\left(\mathbf{X}_{n-2}\right) \right].
\end{equation}
Once this is done, we then solve the problem 
\begin{equation}
\label{eq:Anm1Stars_Ctrl_constr}
A^*_{n-1} 
:=
\underset{
	\footnotesize
	\begin{array}{l}
	\tilde{f}_\mathbf{U}^n \\
	\text{ s.t.:} \\ 
	\constrEq{\mathbf{u}}{n}
	\end{array}
}{\min} \hat{A}(\mathbf{x}_{n-2})
\end{equation}
and obtain  $B^*_{n-1}$ as
\begin{equation}B^*_{n-1}
\begin{array}{ll}
:=
\mathbb{E}_{p^{n-2}_\mathbf{X}}
\left[
A^*_{n-1}
\right].
\end{array}
\end{equation}
To this end we first note that the following identities 
\begin{subequations}
\begin{equation}
\label{eq:Epn_Epnm_Etilfnm}
\mathbb{E}_{p^{n-1}_\mathbf{X}}
\left[ 
\varphi \left(\mathbf{X}_{n-1} \right)
\right]
=
\mathbb{E}_{p^{n-2}_\mathbf{X}}
\left[ 
\mathbb{E}_{\tilde{f}^{n-1}}
\left[ 
\varphi \left(\mathbf{X}_{n-1} \right)
\right]
\right]
\end{equation}
\begin{equation}
\mathbb{E}_{\tilde{f}^{n-1}_\mathbf{X}}
\left[ 
\varphi \left(\mathbf{X}_{n-1} \right)
\right]
=
\mathbb{E}_{\tilde{f}^{n-2}_\mathbf{X}}
\left[ 
\mathbb{E}_{\tilde{f}^{n-1}}
\left[ 
\varphi \left(\mathbf{X}_{n-1} \right)
\right]
\right]
\end{equation}
\end{subequations}
hold for any function $\varphi$ of $\mathbf{X}_{n-1}$. Therefore, by means of
\eqref{eq:ctrl_contrst_min_Bn} and \eqref{eq:Epn_Epnm_Etilfnm}
we obtain, from \eqref{eq:Bnm1_Ctrl_constr}:
\begin{equation}
\label{eq:Opt_Bnm1_shape}
\begin{array}{l}
B_{n-1}   = 
\mathbb{E}_{f^{n-2}}
\left[
\mathcal{D}_{\KL}
\left(
\tilde{f}^{n-1} ||
\tilde{\Ipdf}^{n-1}
\right)
\right]
+
B^\ast_{n} = \\
 =	 
\mathbb{E}_{p^{n-2}_\mathbf{X}}
\left[
\mathcal{D}_{\KL}
\left(
\tilde{f}^{n-1} ||
\tilde{\Ipdf}^{n-1}
\right)
\right]
+
B^\ast_{n} = \\
 =	 
\mathbb{E}_{p^{n-2}_\mathbf{X}}
\left[
\mathcal{D}_{\KL}
\left(
\tilde{f}^{n-1} ||
\tilde{\Ipdf}^{n-1}
\right)
\right]+ \\
\qquad
- 
\mathbb{E}_{p^{n-2}_\mathbf{X}}
\left[
\mathbb{E}_{\tilde{f}^{n-1}}
\left[
\ln \hat{\gamma}\left( \mathbf{X}_{n-1} \right)	
\right]
\right]  
= \\
=	 
\mathbb{E}_{p^{n-2}_\mathbf{x}}
\left[
\underbrace{
\mathcal{D}_{\KL}
\left(
\tilde{f}^{n-1} ||
\tilde{\Ipdf}^{n-1}
\right)
+
\mathbb{E}_{\tilde{f}^{n-1}}
\left[
- \ln \hat{\gamma}\left( \mathbf{X}_{n-1} \right)	
\right]
}_{=:\hat{A}\left(\mathbf{X}_{n-2}\right)  }
\right]  
\end{array}	
\end{equation} 
and the term $\hat{A}\left(\mathbf{x}_{n-2}\right)$ can be recognized.
Now, following the same reasoning we used to compute $\hat{A}(\mathbf{x}_{n-1})$, we explicitly write $\hat{A}(\mathbf{x}_{n-2})$ in compact form as
\begin{equation}
\label{eq:Anm1_final}
\begin{array}{ll}
\hat{A}\left(\mathbf{x}_{n-2}\right) = \\
= 	\mathcal{D}_{\KL}
	\left(
	\tilde{f}^{n-1} ||
	\tilde{\Ipdf}^{n-1}
	\right)
	+
	\mathbb{E}_{\tilde{f}^{n-1}_\mathbf{U}}
	\left[
	\mathbb{E}_{\tilde{f}^{n-1}_\mathbf{X}}
	\left[
	-
	\ln \hat{\gamma}\left( \mathbf{X}_{n-1} \right)	
	\right]
	\right] =  \\
=	\int
	\tilde{f}^{n-1}_\mathbf{U}
	\left\lbrace
	\ln
	\left(
	\frac{
		\tilde{f}^{n-1}_\mathbf{U}
	}{
		\tilde{\Ipdf}^{n-1}_\mathbf{U}
	}
	\right)
	+
	\hat{\omega} \left( \mathbf{u}_{n-1},\,\mathbf{x}_{n-2} \right)
	\right\rbrace
	\, d\mathbf{u}_{n-1}\, ,
\end{array}
\end{equation} 
where $\hat{\omega} \left( \mathbf{u}_{n-1},\,\mathbf{x}_{n-2} \right) = \hat{\alpha} \left( \mathbf{u}_{n-1},\,\mathbf{x}_{n-2} \right)+  \hat{\beta} \left( \mathbf{u}_{n-1},\,\mathbf{x}_{n-2} \right)$ and
\begin{equation}
\begin{array}{rll}
\hat{\alpha} \left( \mathbf{u}_{n-1},\,\mathbf{x}_{n-2} \right)
& :=
\mathcal{D}_{\KL} \left( \tilde{f}^{n-1}_\mathbf{X} || \tilde{\Ipdf}^{n-1}_\mathbf{X} \right)	\\		
\hat{\beta} \left( \mathbf{u}_{n-1},\,\mathbf{x}_{n-2} \right)
& :=
- \mathbb{E}_{\tilde{f}^{n-1}_\mathbf{X}}
\left[
\ln\hat{\gamma}\left( \mathbf{X}_{n-1} \right)	
\right]  	\\
\end{array}
\end{equation}

The last expression for $\hat{A}(\mathbf{x}_{n-2})$ obtained in \eqref{eq:Anm1_final} allows us to use the Lemma~\ref{lem:Constrained_KL} to solve the optimization problem defined in \eqref{eq:Anm1Stars_Ctrl_constr}. Indeed by applying Lemma~\ref{lem:Constrained_KL} with $\mathbf{Z}=\mathbf{U}_{n-1}$, 
$f=\tilde{f}_\mathbf{U}^{n-1}$, 
$g=\tilde{g}_\mathbf{U}^{n-1}$,
$\mathbf{h} = \mathbf{h}_{\mathbf{u},n-1} $,
$\mathbf{H} = \mathbf{H}_{\mathbf{u},n-1} $,
$\hat{\alpha}(\cdot) = \hat{\omega}(\cdot)$,
we get the following solution to the problem in \eqref{eq:Anm1Stars_Ctrl_constr}:
\begin{equation}
\label{eq:opt_ctrl_nm1}
\left( \tilde{f}_\mathbf{U}^{n-1} \right)^\ast = 
\tilde{\Ipdf}_\mathbf{U}^{n-1} 
\frac{e^{-\lbrace 
		\hat{\omega} \left( \mathbf{u}_{n-1},\,\mathbf{x}_{n-2} \right) +
		\langle
		\boldsymbol{\lambda}^*_{\mathbf{u},n-1},
		\mathbf{h}_{\mathbf{u},n-1}\left( \mathbf{u}_{n-1} \right)
		\rangle
		\rbrace}}{e^{1+\lambda^\ast_{\mathbf{u},0,n-1}}}.
\end{equation}
Now, the LM $\lambda^\ast_{\mathbf{u},0,n-1}$ can be obtained by imposing the normalization condition for $\left( \tilde{f}_\mathbf{U}^{n-1} \right)^\ast$. That is,
\begin{equation}
	\begin{split}
		& \exp\{\lambda^\ast_{\mathbf{u},0,n-1} + 1\}  =\\
		& \int \tilde{\Ipdf}_\mathbf{U}^{n-1} e^{-\lbrace \hat{\alpha} \left( \mathbf{u}_{n-1},\,\mathbf{x}_{n-2} \right) +
			\langle
			\boldsymbol{\lambda}^\ast_{\mathbf{u},n-1},
			\mathbf{h}_{\mathbf{u},n-1}\left( \mathbf{u}_{n-1} \right)
			\rangle
			\rbrace}	
		\; d\mathbf{u}_{n-1} = \\
		&  \hat{\gamma}_{\mathbf{u},0,n-1}	
		\left( \mathbf{x}_{n-2} \right).
	\end{split}
\end{equation}
All the other LMs, $\boldsymbol{\lambda}^\ast_{\mathbf{u},n-1}$, can be instead obtained via Lemma \ref{lem:LM_optim}. Moreover, the minimum value for $B_{n-1}^\ast$ corresponding to the above pdf is
\begin{equation}
\label{eq:ctrl_contrst_min_Bnm1}
\begin{array}{rl}
B^\ast_{n-1} 
& = 	 
- \mathbb{E}_{\tilde{f}^{n-2}_\mathbf{X}}
\left[
\sum_{i=0}^{c_{\mathbf{u}}}
\ln 
\hat{\gamma}_{\mathbf{u},i,n-1} \left( \mathbf{X}_{n-2} \right) 	
\right] \\
& = 	 
- \mathbb{E}_{\tilde{f}^{n-2}_\mathbf{X}}
\left[
\ln \hat{\gamma} \left( \mathbf{X}_{n-2} \right)	
\right] .
\end{array}  	
\end{equation}

The proof can then be concluded by observing that at each further backward iteration, the solution 
$\left(\tilde{f}_\mathbf{U}^{k} \right)^*$ has the same shape as 
$\left(\tilde{f}_\mathbf{U}^{n-1} \right)^*$. Indeed, the sub-problems corresponding to each further backward iteration have the exact same structure as the problem solved at the time instant $n-1$. In particular, the problems will have the same structure for the functions $\hat{\alpha}$, $\hat{\beta}$, $\hat{\omega}$, this time evaluated at the previous instants. Moreover, for the last time instant ($n$) the quantity $\hat{\beta} \left( \mathbf{u}_{n},\,\mathbf{x}_{n} \right)$ can be set to $0$ (as there are no constraints at iteration $n+1$) and this is in turn equivalent to have  $\lambda^*_{\mathbf{u},i,n+1} = 0, \forall i$. This completes the proof. 
\qed \\

We are now ready to introduce our  algorithm translating the above theoretical results into a computational tool.
\section{The algorithm}\label{sec:algorithm}
We developed an algorithmic procedure that, by leveraging the technical results introduced above, outputs the solution
$\left\lbrace 	
	\left( \tilde{f}_\mathbf{U}^k \right)^* 
\right\rbrace_{k \in \mathcal{K}}$ 
to Problem \ref{prob:Main_Constr_Ctrl}.
The only inputs that are necessary to the algorithm are $\Ipdf \left( \mathbf{d}_e^n \right)$, extracted from the example dataset and the 
$\tilde f_{\mathbf{X}}^k$'s modeling the plant.
	

\begin{algorithm}[H]
	\caption{Pseudo-code }
	\label{algo:FPD_Analytical_pseudocode}
	\begin{algorithmic}	
			\State \textbf{Inputs:} 
		\State $\Ipdf \left( \mathbf{d}_e^n \right)$ and $\tilde f_{\mathbf{X}}^k$'s
		\State \textbf{Output:}
		\State $\left \lbrace
	\left( \tilde{f}_\mathbf{U}^k \right)^*
	\right  \rbrace_{k \in \mathcal{K}}$ solving Problem \ref{prob:Main_Constr_Ctrl}
		\State \textbf{Initialize} 
		
		\State $\hat{\gamma}_{\mathbf{u},0,n}\left( \mathbf{x}_{\Nh} \right) = 1$
		$\lambda^*_{\mathbf{u},0,n} = 0$, 
		$\lambda^\ast_{\mathbf{u},i,n} = 0$,		
		
		\State
		$\hat{\gamma} \equiv \hat{\gamma}_{\mathbf{u},0,n}$;
		
		\State $\hat{\beta} \left( \mathbf{x}_{\Nh-1},\,\mathbf{u}_{\Nh} \right) = 0$ ;
		\newline
		\For{$ k = \Nh$  to $1$}		 
		\State By backward recursion
		\State
			$\hat{\alpha} \left( \mathbf{u}_{k},	\mathbf{x}_{k-1} \right) 
		 	\gets \int f \left(\mathbf{x}_{k} |	\mathbf{u}_{k},	\mathbf{x}_{k-1} \right) 
		 	\frac{  f \left(\mathbf{x}_{k} |	\mathbf{u}_{k},	\mathbf{x}_{k-1} \right) 
		 	}{ g \left(\mathbf{x}_{k} |	\mathbf{u}_{k},	\mathbf{x}_{k-1} \right) 
		 	} 
		  	\,d\mathbf{x}_{k}
			$
		
		\State
			$ \hat{\beta} \left( \mathbf{u}_{k},	\mathbf{x}_{k-1} \right)
		 	\gets\int 
			f \left(\mathbf{x}_{k} |	\mathbf{u}_{k},	\mathbf{x}_{k-1} \right) 
			\left\lbrace 
				-\ln 
		 		\left(	 \hat{\gamma} \left( \mathbf{x}_{k} \right) 	\right)
			\right\rbrace $
			
		\State
			$ \hat{\omega} \left( \mathbf{u}_{k},	\mathbf{x}_{k-1} \right)
				\gets 	\hat{\alpha} \left( \mathbf{u}_{k},	\mathbf{x}_{k-1} \right) +
					\hat{\beta} \left( \mathbf{u}_{k},	\mathbf{x}_{k-1} \right) $
		\State
			$ \hat{n}_u \left( \mathbf{u}_{k},	\mathbf{x}_{k-1} \right)
		 	\gets g \left( \mathbf{u}_{k} |	\mathbf{x}_{k-1} \right) 
		 	\exp{\left\lbrace
		 	   	- \hat{\omega} \left( \mathbf{u}_{k},	\mathbf{x}_{k-1} \right)
		 	\right\rbrace	}
			$
		\State 
			$ \tilde{\gamma}_0 \left( \mathbf{x}_{k-1} \right)  
		\gets
				\int \hat{n}_u \left( \mathbf{u}_{k},	\mathbf{x}_{k-1} \right)
					\,d\mathbf{u}_{k}  $
		\newline					
			\State
			$ f \left( \mathbf{u}_{k} |	\mathbf{x}_{k-1} \right)
			\gets
			\frac{
				\hat{n}_u \left( \mathbf{u}_{k},	\mathbf{x}_{k-1} \right)
			}{
				\tilde{\gamma}_0 \left( \mathbf{x}_{k-1} \right) 
			} $
			\newline				

			\State Use Lemma \ref{lem:LM_optim} with $\mathcal{Z} := \mathcal{S}(f \left( \mathbf{u}_{k} |	\mathbf{x}_{k-1} \right) )$, $\hat{f}_1 = f$, $\tilde{\mathbf{H}} := \tilde{\mathbf{H}}_{\mathbf{u},k}$,  $\tilde{\mathbf{h}} := \tilde{\mathbf{h}}_{\mathbf{x},k}$, $ \lambda_{0} := \lambda_{\mathbf{u},0,k}$,	$ \boldsymbol{\lambda} := \boldsymbol{\lambda}_{\mathbf{u},k}$, $\tilde{\boldsymbol{\theta}} :=
			\left[ \theta_{0}, \boldsymbol{\theta}^{T}  \right]^{T}
			=  \left[1 +\lambda_{0}, \boldsymbol{\lambda}^{T}  \right]^{T}$ to find the Lagrange multipliers:
			\State $ \boldsymbol{\lambda}^*_{\mathbf{u},k} = \boldsymbol{\lambda}^* \gets  \boldsymbol{\theta}^*$
			\State $ \lambda^*_{\mathbf{u},0,k} 
			\left(
			\mathbf{x}_{k-1}
			\right)=
			\lambda^*_{0} \gets  \theta^*_{0} - 1$ \; 				
			\State	Compute the policy and prepare variables for the next iteration, $k-1$:
			\State $ \left( \tilde{f}_\mathbf{U}^k \right)^* \gets   \frac{f \left( \mathbf{u}_{k} |	\mathbf{x}_{k-1} \right) e^{-
			\langle
			\boldsymbol{\lambda}^*_{\mathbf{u},k}, \mathbf{h}_{\mathbf{u},k}\left( \mathbf{u}_k \right)
			\rangle
			}}{e^{1+\lambda^*_{\mathbf{u},0,k}}}$
			
			\State
			$\hat{\gamma}_{\mathbf{u},i,k}				
			\left(
			\mathbf{x}_{k-1}
			\right) \gets  
			\exp{ \lbrace \lambda^*_{\mathbf{u},i,k} \mathbf{H}_{\mathbf{u},i,k}	\rbrace}  	\quad i = 1,\dots,c_{\mathbf{u}}$
			
			\State 
			$\hat{\gamma}_{\mathbf{u},0,k} = 
			\exp{ \lbrace \theta^*_{0} \rbrace }
			\gets \exp{ \lbrace \lambda^*_{\mathbf{u},0,k} + 1  \rbrace}$ 
			\State
			$ \hat{\gamma}
			\left(
			\mathbf{x}_{k-1}
			\right)
				\gets \exp
				\left[
				\sum_{i=0}^{c_{\mathbf{u}}}
				\ln
				\left(
				\hat{\gamma}_{\mathbf{u},i,k}				
				\left(
				\mathbf{x}_{k}
				\right)
				\right) 	
				\right]
				$
	\EndFor
	\end{algorithmic}
\end{algorithm}

\section{Validation}\label{sec:validation}
\begin{figure}[htbp]
	\centering	
	\includegraphics[width=1\columnwidth]{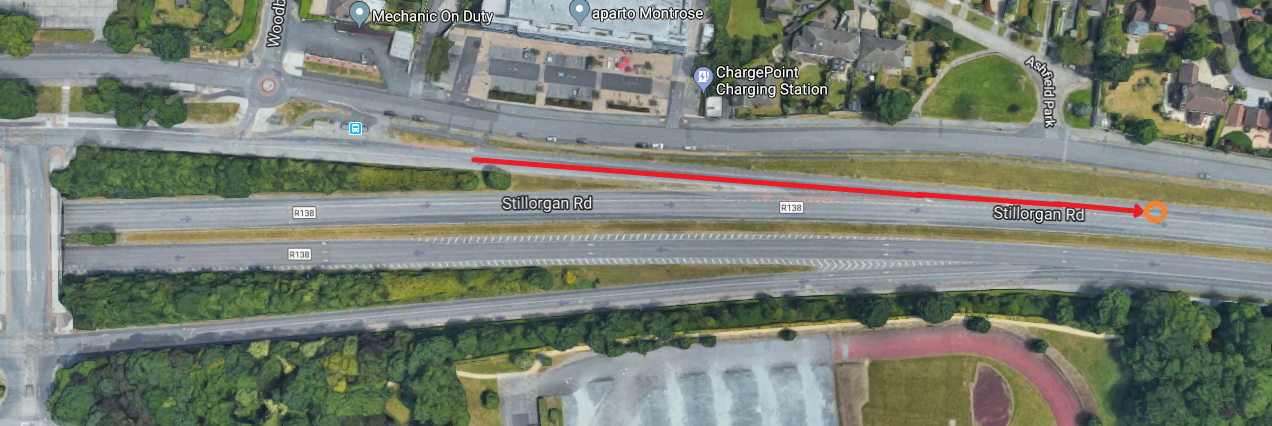}
	\caption{Autonomous driving scenario for Section \ref{sec:validation}: a car that is trying to merge onto a highway. The figure illustrates the stretch of road where the experiments took place. The area is outside the UCD entrance on {\em Stillorgan Road}, Dublin $4$.}
	\label{fig:usecaseAuto}
\end{figure}

We used Algorithm \ref{algo:FPD_Analytical_pseudocode} to synthesize a control policy (from real data) that would allow an autonomous car to merge on a highway.  The scenario considered in our test is described in Fig.\ref{fig:usecaseAuto}. Data were collected using the infrastructure of \cite{griggs2019vehicle}: GPS position, speed, acceleration and jerk were gathered through an OBD2 connection during $100$ test drives.    

\begin{figure}[htbp]
	\centering	
	\includegraphics[width=0.75\columnwidth]{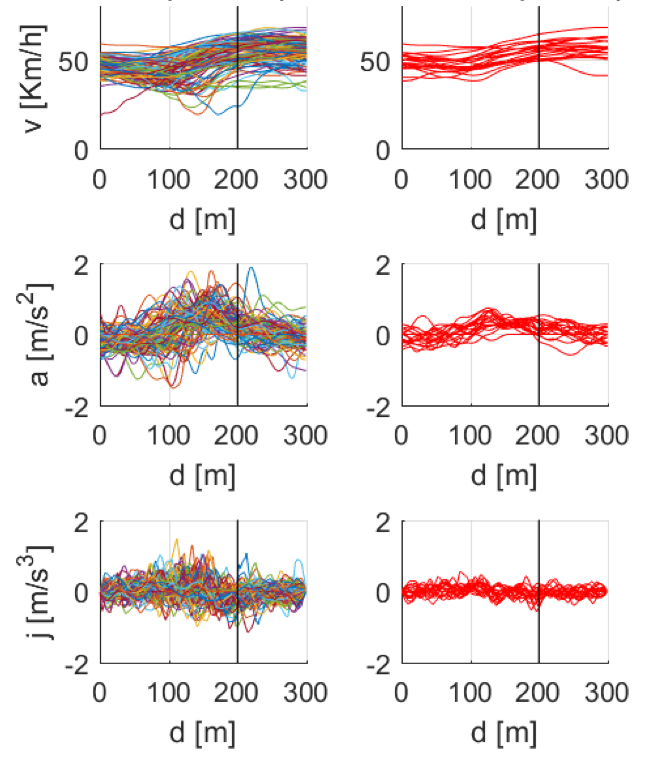}
	\caption{Data collected during the experiments: speed, acceleration, jerk as a function of distance (measured from the beginning of the trip, the UCD entrance). The vertical line in each panel denotes the physical location of the junction highlighted in Fig. \ref{fig:usecaseAuto}. The panels on the left report all the data collected from $100$ trips, while panels on the right report the subset of $20$ trips with the lowest jerk.}
	\label{fig:data}
\end{figure}

The stretch of road we used for our experiments is shown in Fig.\ref{fig:usecaseAuto} and the corresponding data that were collected are in Fig.\ref{fig:data}.

\begin{figure}[htbp]
	\centering	
	\includegraphics[width=1\columnwidth]{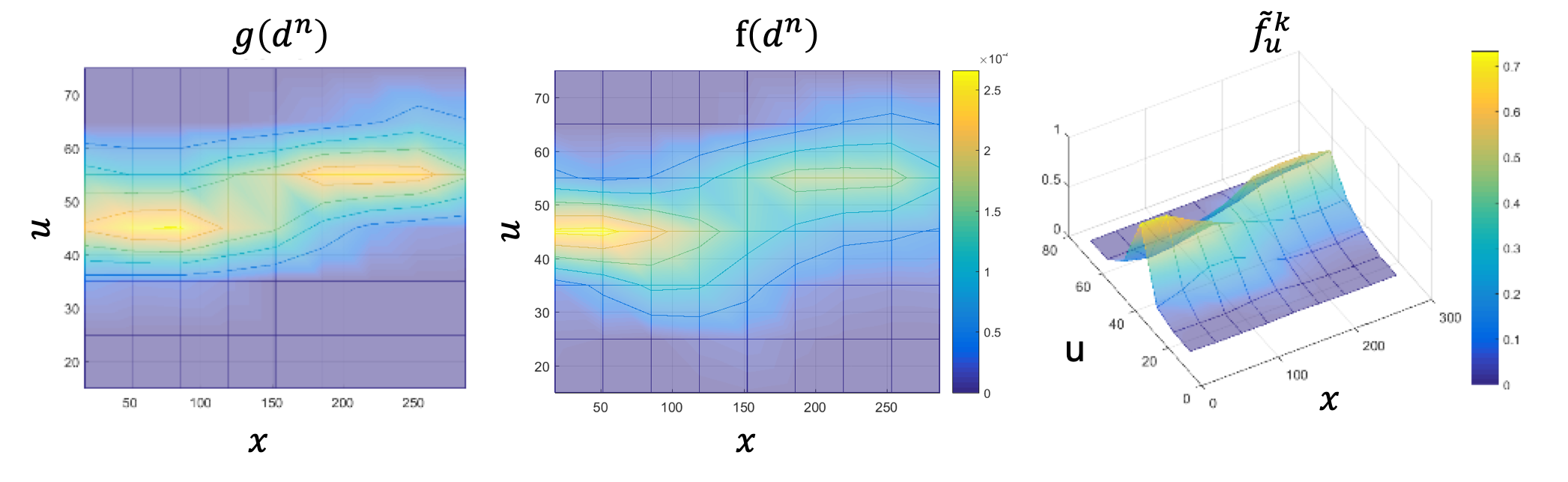}
	\caption{Pdfs extracted from the datasets of Fig. \ref{fig:data}. On the axes, $x$ and $u$ denote the full series of collected distances and speeds.}
	\label{fig:pdfs}
\end{figure}

We used the distance between the the road junction point and the car position as state variable ($\mathbf{x}_k = d(t_k)$) and the car longitudinal speed as control variable ($\mathbf{u}_k = v(t_k)$). From the dataset, we extracted the $20$ trips with the lowest jerk (in red in Fig. \ref{fig:data}). We used this reduced dataset as desired behavior for the car. Given this set-up, we were able to compute both $f(\mathbf{d}^n)$ and $g(\mathbf{d}_e^n)$ from the complete dataset of $100$ trips and the reduced dataset of $20$ trips respectively. These pdfs are shown in Fig. \ref{fig:pdfs}, together with the corresponding control pdf (rightward panel).We also note here that $\text{S}(f(\mathbf{d}^n)) \subseteq \text{S}(g(\mathbf{d}_e^n))$ and this guarantees the absolute continuity of $f(\mathbf{d}^n)$ with respect to $g(\mathbf{d}_e^n)$.

\begin{figure}[htbp]
	\centering	
	\includegraphics[width=1\columnwidth]{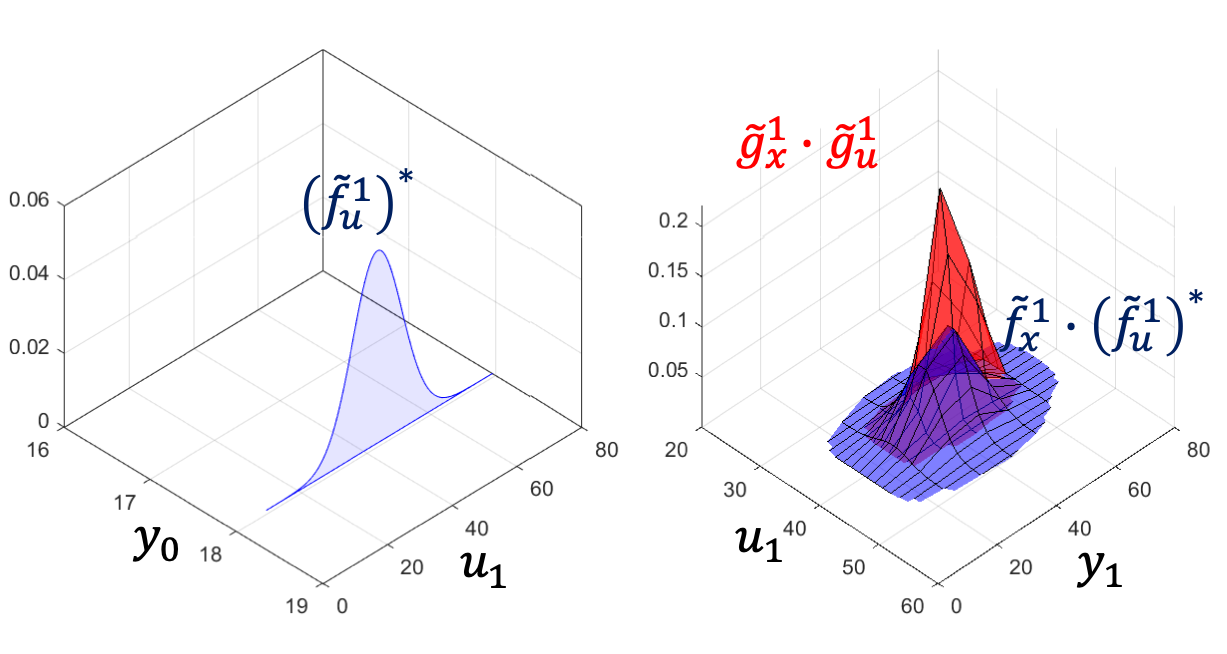}
	\caption{The results obtained using Algorithm \ref{algo:FPD_Analytical_pseudocode}. For the sake of clarity, the results are illustrated at time $k=1$ and are representative of the other time instants. The optimal control pdf (left panel) and the reesulting closed loop pdf (right panel).}
	\label{fig:control}
\end{figure}

Finally, we decided to constraint the variance of the acceleration (the control variable) and solved the resulting Problem \ref{prob:Main_Constr_Ctrl} via Algorithm \ref{algo:FPD_Analytical_pseudocode}. In particular, to make the problem computationally efficient, we approximated all the above pdfs as Gaussian distributions via the Maximum Entropy Principle. Once this was done, we were able to control the closed loop pdf of the system so that it became as close as possible to $g(\mathbf{d}_e^n)$, given the constraint on the variance - see Figure \ref{fig:control}. In the figure, the initial condition was $x_0 = 18$ meters (physically, this is a traffic light outside the UCD gate). Also, the equality constraint was set to have a variance of the closed-loop system higher than the variance of $g(\mathbf{d}_e^n)$ - this is why the closed-loop pdf is flatter.

\section{Conclusions}
We presented an approach to the synthesis of policies from examples. The key technical novelty of the results is the inclusion of actuation constraints in the problem formulation. This in turn yields policies that can be exported to different systems having different actuation capabilities. After presenting the main results we introduced an algorithmic procedure (code is available upon request). If accepted, the presentation will include a sketch of the proofs and a full report of our experimental results, which could not be included here due to space constraints.

\appendix
\section{Sketch of the proofs}
\subsection*{Proof of Property \ref{proper:KLsplit}}
To prove this result we start from the definition of \KL-divergence. In particular:
\begin{equation*}
	\begin{array}{ll}
		\mathcal{D}_{\KL}
		\left( 
		\phi \left( \mathbf{y},\mathbf{z} \right) || 
		g\left( \mathbf{y},\mathbf{z} \right) 
		\right)  := \\
		\quad =
		\int
		\int
		\phi \left( \mathbf{y},\mathbf{z} \right)
		\left[
		\ln \frac{ \phi \left( \mathbf{y},\mathbf{z} \right)}
		{g\left( \mathbf{y},\mathbf{z} \right) } 
		\right]
		\,d\mathbf{y}
		\,d\mathbf{z}	= \\
		\quad =	
		\int
		\int
		\phi \left( \mathbf{z}|\mathbf{y} \right) \,
		\phi \left( \mathbf{y}\right) 
		\left[
		\ln \frac{	\phi \left( \mathbf{z}|\mathbf{y} \right) 
			\phi \left( \mathbf{y}\right)}
		{g \left( \mathbf{z}|\mathbf{y} \right) \,
			g \left( \mathbf{y}\right) } 
		\right]
		\, d\mathbf{y}
		\, d\mathbf{z}	= \\
		\quad = 
		\underbrace{	
			\int
			\int
			\phi \left( \mathbf{z}|\mathbf{y} \right)
			\left[
			\phi \left( \mathbf{y}\right) 
			\ln 
			\frac{
				\phi \left( \mathbf{y}\right) 
			}{
				g \left( \mathbf{y}\right) 
			} 
			\right]
			\, d\mathbf{y}
			\, d\mathbf{z}
		}_{(1)}	\;	+ \\
		\qquad +	
		\underbrace{
			\int 
			\int
			\phi \left( \mathbf{y}\right) 
			\left[
			\phi \left( \mathbf{z}|\mathbf{y} \right)	
			\ln \frac{ 	\phi \left( \mathbf{z}|\mathbf{y} \right)
			}{
				g \left( \mathbf{z}|\mathbf{y} \right)} 
			\right] 
			\, d\mathbf{z}
			\, d\mathbf{y}
		}_{(2)}.	
	\end{array}
\end{equation*}

For the term $(1)$ in the above expression we may continue as follows: 

\begin{equation*}
	\begin{array}{ll}
		&
		\int
		\int
		\phi \left( \mathbf{z}|\mathbf{y} \right)
		\left[
		\phi \left( \mathbf{y}\right) 
		\ln 
		\frac{
			\phi \left( \mathbf{y}\right) 
		}{
			g \left( \mathbf{y}\right) 
		} 
		\right]
		\, d\mathbf{y}
		\, d\mathbf{z}  =\\
		& =
		\int
		\phi \left( \mathbf{z}|\mathbf{y} \right)
		\, d\mathbf{z}
		*
		\left[\int
		\phi \left( \mathbf{y}\right) 
		\ln 
		\frac{
			\phi \left( \mathbf{y}\right) 
		}{
			g \left( \mathbf{y}\right) 
		} 
		d\mathbf{y}\right] 
		\,  =\\
		& = 
		\mathcal{D}_{\KL} \left(\phi \left( \mathbf{y}\right)  || g \left( \mathbf{y}\right)  \right)
	\end{array}
\end{equation*}

where we used Fubini's theorem, the fact that that the term on the first line in square brackets is indepedent on $\mathbf{Z}$ and the fact that $\int
		\phi \left( \mathbf{z}|\mathbf{y} \right)
		\, d\mathbf{z} =1$. 

By using again Fubini's theorem, for the term $(2)$ instead we have:
\begin{equation*}
	\begin{array}{ll}
		& 
		\int 
		\int
		\phi \left( \mathbf{y}\right) 
		\left[
		\phi \left( \mathbf{z}|\mathbf{y} \right)	
		\ln \frac{ 	\phi \left( \mathbf{z}|\mathbf{y} \right)
		}{
			g \left( \mathbf{z}|\mathbf{y} \right)} 
		\right] 
		\, d\mathbf{z}
		\, d\mathbf{y}
		=	\\
		& =
		\int
		\phi \left( \mathbf{y}\right) 
		\left[ 
		\int
		\phi \left( \mathbf{z}|\mathbf{y} \right)	
		\ln \frac{ 	\phi \left( \mathbf{z}|\mathbf{y} \right)
		}{
			g \left( \mathbf{z}|\mathbf{y} \right)} 
		\, d\mathbf{z}
		\right]
		\, d\mathbf{y}=  \\
		& =
		\int
		\phi \left( \mathbf{y}\right) 
		\left[	
		\mathcal{D}_{\KL} \left( 
		\phi \left( \mathbf{z}|\mathbf{y} \right)|| 
		g \left( \mathbf{z}|\mathbf{y} \right)\right) 
		\right]
		\, d\mathbf{y}= \\
		& =
		\mathbb{E}_{ \phi ( \mathbf{Y} ) }
		\left[	
		\mathcal{D}_{\KL} \left( 	\phi(\mathbf{z}|\mathbf{Y}) || \, g(\mathbf{z}|\mathbf{Y}) \right)
		\right],
	\end{array}
\end{equation*}
thus proving the result.
\qed

\subsection*{Proof of Lemma \ref{lem:Constrained_KL}}
We prove the result in two steps. First, we rewrite the cost function $\mathcal{L}(f)$ and consider the corresponding augmented Lagrangian. Then, we make use of the Euler-Lagrange (EL)  stationary conditions to find $f^\ast_{\mathbf{Z}}(\mathbf{z})$ (in what follows we omit the dependencies of functions and pdfs on the random variable $\mathbf{z}$ whenever this is clear from the context). 

As a first step, note that the cost function $\mathcal{L}(f)$ of the constrained optimization problem in (\ref{eqn:probl_Lemma_1}) can be conveniently re-written as $\mathcal{L}(f) = \int f\left[\ln\left(\frac{f}{\Ipdf}\right)+\alpha\right] \,d\mathbf{z}$. 
Then, the augmented Lagrangian takes the following form:
\begin{equation*}
\begin{split}
&\mathcal{L}_{aug}
\left( f,\lambda_0, \boldsymbol{\lambda} \right) := \\
&
	\int f \; 
	\left[
	\ln
	\left(
	\frac{f}{
		\Ipdf}
	\right)
	+ 
	\alpha
	\right] \,
	d\mathbf{z} + \lambda_0
	\left(
	\int f \; 
	\mathds{1}_{\mathcal{S}  \left( \mathbf{Z} \right) }
	\; d\mathbf{z} - 1
	\right) \\
	& +
	\langle 
	\boldsymbol{\lambda},
	\int f \; 
	\mathbf{h}\left( \mathbf{z} \right)\; d\mathbf{z} 
	-\mathbf{H}
	\rangle,
	\end{split}
	\end{equation*}
where $\lambda_0$ and $\boldsymbol{\lambda}:= [\lambda_1,\ldots,\lambda_{c_\mathbf{z}}]^T$ are the (non-negative) Lagrange multipliers (LMs) corresponding to the constraints of the optimization problem. In turn, the above expression can be re-written as
\begin{equation}
\label{eqn:lagrangian_Lema_1}
\begin{split}
&\mathcal{L} _{aug}
\left( f,\lambda_0, \boldsymbol{\lambda} \right)
 =
	\int f \; 
	\left[
	\ln
	\left(
	\frac{f}{
		\Ipdf}
	\right)
	+ 
	\alpha
	+
	\lambda_0 
	+ 
	\langle
	\boldsymbol{\lambda},
	\mathbf{h}\left( \mathbf{z} \right)
	\rangle
	\right] \,
	d\mathbf{z}	\\
	&- \lambda_0 
	- \langle\boldsymbol{\lambda},\mathbf{H}\rangle
	\end{split}
	\end{equation}
Now, we let
\begin{equation}\label{eqn:alpha_expression}
\tilde{\alpha} (\mathbf{z})
=
\alpha (\mathbf{z} )
+
\lambda_0 
+ 
\langle
\boldsymbol{\lambda},
\mathbf{h}\left( \mathbf{z} \right)
\rangle
\end{equation}
and  make use of the EL stationary conditions to find the optimal solution. First, we consider the EL stationary condition with respect to the pdf $f$. These conditions can be written in terms of the quantity under the integral in (\ref{eqn:lagrangian_Lema_1}), i.e. in terms of $ l(f) := f \; \left[ \ln \left( \frac{f}{g} \right) + \tilde{\alpha} \right] = f \; \left[ \ln \left(f \right) - \ln \left(g \right) + \tilde{\alpha} \right] $. In particular, by imposing the stationary condition we obtain:
\begin{equation}
\frac{\partial l(f)}{\partial f }  = \ln \left( \frac{f}{g} \right) + \tilde{\alpha} + 1 = 0.
\end{equation}
Therefore, it follows that all the optimal solution candidates must be of the form: 
\begin{equation}
	f(\mathbf{z} ) = g(\mathbf{z} )e^{-\lbrace 1 + \tilde{\alpha}(\mathbf{z} )\rbrace},
\end{equation}
which, by definition of $\tilde{\alpha}$, becomes
\begin{equation}\label{eqn:candidate}
f(\mathbf{z}) = 
g(\mathbf{z} )\frac{e^{-\lbrace \alpha(\mathbf{z}) + \langle \boldsymbol{\lambda},\mathbf{h}\left( \mathbf{z} \right) \rangle \rbrace}}
		{e^{1+\lambda_0}}.
\end{equation}
Note that the above candidates are a function of the LMs. These can be computed by applying the EL stationary condition with respect to $\lambda_0,\lambda_1,\ldots,\lambda_{c_\mathbf{z}}$. This yields the following set of additional conditions:
\begin{equation}\label{eqn:lambdas_conditions}
	\frac{\partial\mathcal{L} _{aug}
		\left(
		f,
		\lambda_0,
		\boldsymbol{\lambda}
		\right)}{\partial \lambda_i}
		= 0, \ \ \ i =0,\ldots, c_{\mathbf{z}}.
\end{equation}
That is, (\ref{eqn:lambdas_conditions}) imply that the LMs  associated to the constraints must satisfy:
\begin{equation}
\int 
g(\mathbf{z} )\frac{e^{-\lbrace \alpha(\mathbf{z}) + \langle \boldsymbol{\lambda},\mathbf{h}\left( \mathbf{z} \right) \rangle \rbrace}}
		{e^{1+\lambda_0}} 	\tilde{\mathbf{h}}_i\left( \mathbf{z} \right)\; d\mathbf{z} 
= 
\tilde{\mathbf{H}}_i
\quad 
i = 0,\dots,c_{\mathbf{z}},
\end{equation}
which was obtained by replacing the expression of the optimal solution candidate (\ref{eqn:candidate}) in (\ref{eqn:lambdas_conditions}).

Now, the above set of equations can be solved via Lemma \ref{lem:LM_optim} and here we let $\lambda_0^\ast$, $\boldsymbol{\lambda}^\ast$ be the resulting values of LMs. By substituting the optimal LMs into the expression of the optimal solution candidates yields:
\begin{equation*}
f^\ast:=f^\ast(\mathbf{z}) = 
g(\mathbf{z} )\frac{e^{-\lbrace \alpha(\mathbf{z}) + \langle \boldsymbol{\lambda}^\ast,\mathbf{h}\left( \mathbf{z} \right) \rangle \rbrace}}
		{e^{1+\lambda_0^\ast}}.
\end{equation*}
The proof is then concluded by noticing that $f^\ast(\mathbf{z})$ is indeed the optimal solution since the Lagrangian is convex in $f$. To show convexity, it suffices to consider the second derivative of $l(f)$ and to observe that this is always positive definite (indeed
$\frac{\partial^2 l}{\partial f^2 } = \frac{\partial \left( \ln \left( f \right) + \tilde{\alpha} + 1  \right)}{\partial f }  = \frac{1}{f}>0$). 

Finally, the second part of the result follows from evaluating $\mathcal{L} 
\left( f^\ast \right)$. Indeed:
\begin{equation}
\begin{array}{rl}
\mathcal{L} 
\left( f^\ast \right) 
{\tiny }	& = \int f^\ast 
	\left[
	\ln \frac{ \Ipdf e^{-\lbrace \alpha +	
			\langle \boldsymbol{\lambda}^\ast, \mathbf{h} \rangle \rbrace}}{e^{1+\lambda_0^\ast} \Ipdf} 
	+ \alpha 
	\right]
	\,d\mathbf{z} = \\
	& = -\int f^* 
	\left( 1+\lambda_0^\ast +	\langle\boldsymbol{\lambda}^\ast,\mathbf{h} \rangle	
	\right)
	\,d\mathbf{z} =\\
	& = -\left( 1+\lambda_0^\ast +\langle \boldsymbol{\lambda}^\ast,\mathbf{H}	\rangle
	\right),	
\end{array}
\end{equation} 
and this completes the proof.\qed
\subsection*{Proof of  Lemma \ref{lem:LM_optim}}
We prove this result by showing that: 
(i) $\mathcal{J}
\left( \tilde{\boldsymbol{\theta}} \right)$ is strictly convex; (ii) its minimizer must satisfy the set of equations (\ref{eq:constr_tilde_Theta}). 

The proof of statement (ii) comes directly from the evaluation of the first order stationary condition. Indeed, any optimal candidate, say $\tilde{\boldsymbol{\theta}}^\ast$, must satisfy $\nabla \mathcal{J}  \left( \tilde{\boldsymbol{\theta}^\ast}\right) =0$.  Now, computing $\nabla \mathcal{J} \left( \tilde{\boldsymbol{\theta}^\ast} \right)$ yields 
\begin{equation}\label{eqn:equiv_condition}
\begin{array}{l}
\tilde{\mathbf{H}}-
	\int_{\mathcal{Z}} 
	\hat{f}_1 \left( \mathbf{z} \right) 	
	e^{- 
		\langle \tilde{\boldsymbol{\theta}^\ast}  
		,\; 	\tilde{\mathbf{h}} \left( \mathbf{z} \right)
		\rangle
	}	
	\tilde{\mathbf{h}} \left( \mathbf{z} \right)
	\, d\mathbf{z} =\\
	 \tilde{\mathbf{H}}-\int_{\mathcal{Z}} 
	\hat{f}_2\left(\mathbf{z},\tilde{\boldsymbol{\theta}^\ast}\right)
	\tilde{\mathbf{h}} \left( \mathbf{z} \right)
	\, d\mathbf{z}= \mathbf{0},\\	
\end{array}
\end{equation}
where we used the definition of $\hat{f}_2$ to obtain the second equality.
That is, (\ref{eqn:equiv_condition}) immediately implies that any candidate minimizer of the optimization problem in (\ref{eq:Index_Inique}) must fulfil the set of equations (\ref{eq:constr_tilde_Theta}). 

In order to prove strict convexity (i.e. statement (i)) we compute the Hessian of
$\mathcal{J} \left( \tilde{\boldsymbol{\theta}} \right)$ and show that this is strictly positive definite in $\tilde{{\Theta}}$. Indeed, computing the Hessian yields
\begin{equation}
\begin{array}{l}
\left[\nabla^2
\mathcal{J} 
\left( \tilde{\boldsymbol{\theta}}  \right)
\right]
=
\int_{\mathcal{Z}} 
\left[
\tilde{\mathbf{h}} \left( \mathbf{z} \right)
\otimes \tilde{\mathbf{h}} \left( \mathbf{z} \right)
\right]
\hat{f}_1 \left( \mathbf{z} \right) 	
	e^{- 
		\langle \tilde{\boldsymbol{\theta}}  
		,\; 	\tilde{\mathbf{h}} \left( \mathbf{z} \right)
		\rangle
	}	
\, d\mathbf{z},
\end{array}
\end{equation}
where $\otimes$ denotes the external product between tensors.  
Now, since the equations in (\ref{eq:constr_tilde_Theta}) are 
algebraically independent, we have (see Definition \ref{def:algebrIndep}):
\begin{equation}
\begin{array}{l}
\exists S \subset \mathcal{Z}: \forall 
\tilde{\mathbf{v}}
\in \mathbb{R}^{c_{\mathbf{z}}}-\{\mathbf{0}\} \\
\langle
\left[\nabla^2
\mathcal{J} 
\left(  \tilde{\boldsymbol{\theta}}	\right)
\right]\tilde{\mathbf{v}}, \,
\tilde{\mathbf{v}}
\rangle
=\\
\int_S 
\langle 
\tilde{\mathbf{v}},
\,\tilde{\mathbf{h}} \left( \mathbf{z} \right)
\left( \mathbf{z} \right) 
\rangle^2
\hat{f}_1 \left( \mathbf{z} \right) 	
e^{- 
	\langle \tilde{\boldsymbol{\theta}}  
	,\; 	\tilde{\mathbf{h}} \left( \mathbf{z} \right)
	\rangle
}	
\, d\mathbf{z}
> 0.\\
\end{array}
\end{equation}	
This implies that $\tilde{\boldsymbol{\theta}}^\ast$ is the unique minimizer of the optimization problem, thus concluding the proof. \qed


\end{document}